\documentclass{daj}

\usepackage{amsmath,amscd,amsfonts,amssymb}
\usepackage{mathrsfs,dsfont}
\usepackage{graphicx}
\usepackage{color}
\usepackage{mathtools}
\usepackage{tikz-cd}

\usepackage{tikz}
\usetikzlibrary{decorations.pathreplacing}

\numberwithin{equation}{section}
\numberwithin{figure}{section}

\usepackage{amsthm,amsmath,amssymb,mathtools}

\usepackage{graphicx}

\usepackage{todonotes}

\newcommand\Z{\mathbb{Z}}

\theoremstyle{plain}
\newtheorem{theorem}{Theorem}
\newtheorem{lemma}[theorem]{Lemma}

\newtheorem{proposition}[theorem]{Proposition}

\theoremstyle{definition}
\newtheorem{definition}[theorem]{Definition}

\theoremstyle{remark}
\newtheorem{remark}[theorem]{Remark}

\dajAUTHORdetails{%
  title = {Additive Energies on Discrete Cubes}, 
  author = {Jaume de Dios Pont, Rachel Greenfeld, Paata Ivanisvili and Jos\'e Madrid},
  plaintextauthor = {Jaume de Dios Pont, Rachel Greenfeld, Paata Ivanisvili, Jos\'e Madrid},
  plaintexttitle = {Additive Energies on Discrete Cubes}, 
  runningauthor = {J. de Dios Pont, R. Greenfeld, P.Ivanisvili and J. Madrid},
  copyrightauthor = {J. de Dios Pont, R. Greenfeld, P.Ivanisvili and J. Madrid},
 }  

\dajEDITORdetails{%
   year={2023},
   number={13},
   received={19 December 2021},   
   revised={18 January 2023},    
   published={1 September 2023},  
   doi={10.19086/da.84737},       
}   

\begin{document}

\begin{frontmatter}[classification=text]

\title{Additive Energies on Discrete Cubes} 


\author[J. de Dios]{Jaume de Dios Pont}

\author[R. Greenfeld]{Rachel Greenfeld\thanks{Partially supported by the Eric and Wendy Schmidt Postdoctoral Award, the  AMIAS Membership  and NSF  grants  DMS-2242871 and   DMS-1926686.}}

\author[P. Ivanisvili]{Paata Ivanisvili\thanks{Partially supported by the NSF grants DMS-2152346 and CAREER-DMS-215240.}}

\author[J. Madrid]{Jos\'e Madrid}

\begin{abstract}
We prove that for $d\geq 0$ and $k\geq 2$, for any subset $A$ of a discrete cube $\{0,1\}^d$, the $k-$higher energy of $A$ (i.e., the number of $2k-$tuples $(a_1,a_2,\dots,a_{2k})$ in $A^{2k}$ with $a_1-a_2=a_3-a_4=\dots=a_{2k-1}-a_{2k}$) is at most $|A|^{\log_{2}(2^k+2)}$, and $\log_{2}(2^k+2)$ is the best possible exponent. We also show that if $d\geq 0$ and $2\leq k\leq 10$, for any subset $A$ of a discrete cube $\{0,1\}^d$, the $k-$additive energy of $A$ (i.e., the number of $2k-$tuples $(a_1,a_2,\dots,a_{2k})$ in $A^{2k}$ with $a_1+a_2+\dots+a_k=a_{k+1}+a_{k+2}+\dots+a_{2k}$) is at most $|A|^{\log_2{ \binom{2k}{k}}}$, and $\log_2{ \binom{2k}{k}}$ is the best possible exponent.  We discuss the analogous problems for the  sets $\{0,1,\dots,n\}^d$ for $n\geq2$.
\end{abstract}
\end{frontmatter}


\section{Introduction}

The {additive energy} $E(A)$ of a finite subset $A$ of an additive group $G$ is defined as the number of quadruples $(a_1,a_2,a_3,a_4) \in A \times A \times A \times A$ such that $a_1+a_2=a_3+a_4$ (see \cite{TaoVu}). Observe that for any triple $(a_1,a_2,a_3)$ there is at most one $a_4$ such that $a_1+a_2=a_3+a_4$, so we have the trivial upper bound $E(A) \leq |A|^3$ (here $|A|$ denotes the cardinality of $A$). This bound is attained, for example, when $A$ is itself a finite group. Considering the diagonal solutions $a_1=a_3$ and $a_2=a_4$ we also observe the trivial lower bound $E(A)\geq |A|^2$.

\subsection{Higher energies}
We define the $k-$higher energy of a set $A\subseteq\{0,1\}^d\subset \mathbb Z^d$ by
$$
\widetilde E_k(A):=|\{(a_1,a_2,\dots,a_{2k-1},a_{2k})\in A^{2k}:a_1-a_2=a_3-a_4=\dots=a_{2k-1}-a_{2k}\}|.
$$
This has been studied by many authors, see \cite{Sh}, \cite{ScSh}. In this case we have the trivial bounds $|A|^k\leq \widetilde E_k(A)\leq |A|^{k+1}$.

\begin{theorem}\label{main thm 2}  
    Let $d \geq 0$, $k\geq 2$, and let $A \subset \{0,1\}^d$.  
    Then $\widetilde E_{k}(A) \leq |A|^{q_k}$, where $q_{k} \coloneqq \log_2{(2^k+2)}$.  
    Furthermore, the exponent $q_k$ cannot be replaced by any smaller quantity.
\end{theorem}
\begin{remark}
This Theorem extends a result obtained by Kane--Tao \cite[Theorem 7]{KT} for $k=2$.
\end{remark}
The second claim in our Theorem \ref{main thm 2} follows considering the case $A = \{0,1\}^d$, 
in this case we have $|A|=|\{0,1\}|^d = 2^d$ and $\widetilde E_{k}(\{0,1\}^d)=(2^k+2)^d$.

\subsection{k-additive energies}
\label{sub:dummy}
We discuss another generalization of Kane--Tao result \cite[Theorem 7]{KT}. 
We define the $k$-additive energy $E_{k}(A)$ of a subset $A$ of an additive group $G$ as the number of $2k-$tuples $(a_1,a_2,\dots,a_{2k})$ in $A^{2k}$ with $a_1+a_2+\dots+a_k=a_{k+1}+a_{k+2}+\dots+a_{2k}$. In this case the trivial bounds are  $|A|^k\leq E_k(A)\leq |A|^{2k-1}$, and we have the following refinement in the cube $\{0,1\}^{d}.$

\begin{theorem}\label{main thm 1}  
    Let $d \geq 0$, $2\leq k\leq 10$, and let $A \subset \{0,1\}^d$.  
    Then $E_{k}(A) \leq |A|^{p_k}$, where $p_{k} \coloneqq \log_2 {\binom{2k}{k}}$.  
    Furthermore, the exponent $p_k$ cannot be replaced by any smaller quantity.
\end{theorem}

\begin{remark}
Theorem \ref{main thm 1} also extends a result obtained by Kane--Tao (\cite[Theorem 7]{KT}). 
\end{remark}

From the well-known bounds for the central binomial coefficient
$\frac{4^k}{2\sqrt{\pi k}}\leq {\binom{2k}{k}}\leq \frac{4^k}{\sqrt{\pi k}}$, one recovers
\begin{equation}\label{eq: upper bound for p}
p_k<2k-1.
\end{equation}

As previously, the second claim in our Theorem \ref{main thm 1} follows considering the case $A = \{0,1\}^d$, since in this case we have $|A|=|\{0,1\}|^d = 2^d$ and $E_{k}(A) =\left[\sum_{i=0}^{k} {\binom{k}{i}}^{2}\right]^d={\binom{2k}{k}}^d$.  
We prove this theorem by induction on $d$ together with the following subtle inequality for Legendre polynomials.

\begin{lemma}
    \label{elem}
     Let $2\leq k\leq 10$ and $p_k=\log_2 {\binom{2k}{k}}$.  
     If $a,b \geq 0$, then

\begin{align}
    \label{elem ineq}
\sum_{j=0}^k \binom{k}{j}^{2} a^{p_k\frac{k-j}{k}} b^{p_k \frac j  k} \le (a+b)^{p_k}.
\end{align}
\end{lemma}

The polynomials  $Q_{k}(t)$, $k\geq 0$, defined by

\begin{align*}
Q_{k}(t)=\frac{1}{2^{k} k!} \frac{\mathrm{d}^{k}}{\mathrm{d}t^{k}}(t^{2}-1)^{k}=\frac{1}{2^{k}} \sum_{j=0}^{k} \binom{k}{j}^{2} (t-1)^{k-j}(t+1)^{j}   
\end{align*}
are called Legendre polynomials. They are orthogonal with respect to Lebesgue measure on the interval $[-1,1]$, each $Q_{k}(t)$ has degree $k$, and they satisfy normalization constraint $Q_{k}(1)=1$. Dividing both sides of (\ref{elem ineq}) by $a^{p_{k}}$ (without loss of generality assume $a\neq 0$), then (\ref{elem ineq}) takes the  form $(y-1)^{k}Q_{k}\left(\frac{y+1}{y-1}\right)\leq (1+y^{k/p_{k}})^{p_{k}}$ with $y=(b/a)^{p_{k}/k}\geq 0$. If we let $t :=\frac{y+1}{y-1}$ (without loss of generality assume  $y\geq 1$), then (\ref{elem ineq}) is the same as 
\begin{align*}
    Q_{k}(t) \leq \left(\left(\frac{t-1}{2}\right)^{\frac{k}{p_{k}}} +\left(\frac{t+1}{2}\right)^{\frac{k}{p_{k}}} \right)^{p_{k}} \quad \text{for all} \quad t\geq 1. 
\end{align*}
This explains the reason  we call Lemma~\ref{elem}
 the inequality for Legendre polynomials.

\subsection{More general discrete cubes}
Let $d\geq0$. Let us consider additive energies of subsets of general discrete cubes\footnote{A related problem about the lower bound for the size of sumsets of subsets of the general discrete cube was studied, e.g., in \cite[Theorem 5]{Bour}. }  $\{0,1,\dots,n\}^d$. Let $t_n$ be the smallest number such that $$E_2(A) \le |A|^{t_n}$$ for all $A\subseteq \{0,1,\dots,n\}^d$. We have seen that in both Theorem \ref{main thm 2} and Theorem \ref{main thm 1} 
we have $q_k=\frac{\log \widetilde E_{k}(\{0,1\}^d)}{\log|\{0,1\}|^d}$ and $p_k=\frac{\log E_{k}(\{0,1\}^d)}{\log|\{0,1\}|^d}$. 
Thus, one could  a-priori  expect a  similar phenomenon for the additive energy of $\{0,1,\dots,n\}^d$. 
However, it turns out that this is not the case in general, not even for the discrete cube $\{0,1,2\}^d$.

\begin{proposition}\label{Thm: 0,1,2}
The following inequality holds
$$
t_2>\frac{\log E_2(\{0,1,2\}^d)}{\log |\{0,1,2\}^d|}.
$$
\end{proposition}
Although finding the precise values of the optimal powers $t_n$ for general discrete cubes $\{0,1,\dots,n\}^d$ seems to be a difficult problem, 
we obtain some bounds describing the asymptotic behavior of $t_n$ as $n$ goes to infinity.

\begin{proposition}\label{Thm larger cubes}
If $n=2m-1$, then 
$$3\geq t_n\geq \log_{2m}\left(\frac{16m^3+2m}{3}\right)>3-\frac{\log(3/2)}{\log(2m)}.$$
If $n=2m$, then
$$3\geq t_n\geq \log_{2m}\left(\frac{16m^3+24m^2+14m+3}{3}\right)>3-\frac{\log(3/2)}{\log(2m)}.$$
\end{proposition}

\section{Proof of Theorem \ref{main thm 2}}

\newcommand{\acorr}{\star}

The proof of Theorem \ref{main thm 2} proceeds via induction on $d$. Observe that the result is trivial for $d=0$.  Assume now that $d \geq 1$ and that the result has been established for $d-1$. Any set $A \subseteq \{0,1\}^d$ can be written as
$$ A = (A_0 \times \{0\}) \uplus (A_1 \times \{1\})$$
for some $A_0,A_1 \subseteq \{0,1\}^{d-1}$, where $\uplus$ means disjoint union. Then we have
\begin{align}\label{key formula 2}
\widetilde E_k(A)&= 
|\{(a_1,a_2,\dots,a_{2k})\in (A_0\times A_1)^{k}: a_1-a_2=a_3-a_4=\dots=a_{2k-1}-a_{2k}\}|\nonumber\\
&\ \ \ +|\{(a_1,a_2,\dots,a_{2k})\in (A_1\times A_0)^{k}: a_1-a_2=a_3-a_4=\dots=a_{2k-1}-a_{2k}\}|\nonumber\\
&\ \ \ +\sum_{i=0}^{k} {\binom{k}{i}}|\{ {(a_1,a_2,\dots,a_{2k})} \in (A^2_0)^i\times (A^{2}_1)^{k-i}\nonumber\\
&\ \ \ \ \ \ \ \ \ \ \ \ \ \ \ \ \ \ \ \ \ \ : a_1-a_2=a_3-a_4=\dots=a_{2k-1}-a_{2k}\}|\nonumber\\
&=:C_1+C_2+\widetilde E_k(A_0)+\widetilde E_k(A_1)+ \sum_{i=1}^{k-1} {\binom{k}{i}}C_{i,k}.
\end{align}
The next proposition plays a fundamental role in our proof.
\begin{proposition}\label{bounds for c i k 2}
For all $1\leq i\leq k-1$ we have that
\begin{equation*}
C_{i,k}\leq |A_0|^{\frac{i}{k}q_k}|A_{1}|^{\frac{k-i}{k}q_k}.
\end{equation*}
Moreover, we have that
$$
C_1\leq |A_0|^{\frac{q_k}{2}}|A_1|^{\frac{q_k}{2}}  \ \text{and}\ \ C_2\leq |A_0|^{\frac{q_k}{2}}|A_1|^{\frac{q_k}{2}}.
$$
\end{proposition}
\begin{proof}[Proof of Proposition \ref{bounds for c i k 2}]
We observe that
$$
\widetilde E_k(A):=\sum_{x\in \Z^d}(\chi_{A}\acorr \chi_{A})^k(x),
$$
where $\chi_A$ denotes the characteristic function of the set $A$, and $f\acorr g$ denotes the correlation of the functions $f$ and $g$ defined by $f\acorr g(x):=\sum_{y\in\Z^d}f(y)g(x+y)$ \cite[Equation 7]{ScSh}. Moreover, by H\"older's inequality we have 
\begin{align*}
    C_{i,k}&=\sum_{x\in\Z^d}(\chi_{A_0}\acorr\chi_{A_0})^i(x)(\chi_{A_1}\acorr\chi_{A_1})^{k-i}(x)\\
    &\leq  \left(\sum_{x\in\Z^d}(\chi_{A_0}\acorr\chi_{A_0})^k(x)\right)^{\frac{i}{k}}\left(\sum_{x\in\Z^d}(\chi_{A_1}\acorr\chi_{A_1})^k(x)\right)^{\frac{k-i}{k}}\\
    &=\widetilde E^{\frac{i}{k}}_k(A_0)\widetilde E^{\frac{k-i}{k}}_k(A_1)\\
    &\leq |A_0|^{\frac{q_ki}{k}}|A_1|^{\frac{q_k(k-i)}{k}}.
\end{align*}
The first identity follows from the facts that
$\chi_{A_0}\acorr\chi_{A_0}(x)$ counts the number of pairs $(y,z)\in A^2_0$ such that $z-y=x$, and $\chi_{A_1}\acorr\chi_{A_1}(x)$ counts the number of pairs $(y,z)\in A^2_1$ such that $z-y=x$. We define
$$ f\bullet g :=\sum_{\substack{a_1,a_2,\dots,a_k\in \{0,1\}^d\\ b_1,b_2,\dots,b_k\in\{0,1\}^d\\ a_1-b_1=a_2-b_2=\dots=a_k-b_k}}f(a_1)f(a_2)\dots f(a_k)g(b_1)g(b_2)\dots g(b_k).
$$
Then
\begin{align*}
&f\bullet g=\\
&\sum_{c_2,c_3,\dots,c_k\in\{-1,0,1\}^d}\left(\sum_{\substack{a_1\in\{0,1\}^d\\  a_1+c_i\in\{0,1\}^d}}f(a_1)f(a_1+c_2)f(a_1+c_3)\dots f(a_1+c_k)\right)\\
&\ \ \ \times\left(\sum_{\substack{b_1\in\{0,1\}^d\\b_1+c_i\in\{0,1\}^d}}g(b_1)g(b_1+c_2)g(b_1+c_3)\dots g(b_1+c_k)\right).
\end{align*}
Therefore, by the Cauchy-Schwarz inequality we obtain
\begin{align*}
    C_1= \chi_{A_0}\bullet\chi_{A_1}&\leq  (\chi_{A_0}\bullet\chi_{A_0})^{1/2}( \chi_{A_1}\bullet\chi_{A_1})^{1/2}\\
    &=\widetilde E^{1/2}_k(A_0)\widetilde E^{1/2}_{k}(A_1)\leq |A_0|^{\frac{q_k}{2}}|A_1|^{\frac{q_k}{2}}. 
\end{align*}
Similarly $C_2 \leq |A_0|^{\frac{q_k}{2}}|A_1|^{\frac{q_k}{2}}$.
\end{proof}

Then, from \eqref{key formula 2}, using Proposition \ref{bounds for c i k 2}
we obtain
\begin{align*}
    \widetilde E_k(A)&=C_1+C_2+\widetilde E_k(A_0)+\widetilde E_k(A_1)+ \sum_{i=1}^{k-1} {\binom{k}{i}}C_{i,k}\\
    &\leq 2|A_0|^{\frac{q_k}{2}}|A_1|^{\frac{q_k}{2}}+ \sum_{i=0}^{k} {\binom{k}{i}}|A_0|^{\frac{i}{k}q_k}|A_{1}|^{\frac{k-i}{k}q_k}\\
    &=2|A_0|^{\frac{q_k}{2}}|A_1|^{\frac{q_k}{2}}+(|A_0|^{\frac{q_k}{k}}+|A_1|^{\frac{q_k}{k}})^k. 
\end{align*}
Thus, to complete the inductive argument, it is enough to prove that for $x=|A_0|$ and $y=|A_1|$ one has
\begin{equation}\label{Key elem ineq for higher energy}
2x^{\frac{q_k}{2}}y^{\frac{q_k}{2}}+(x^{\frac{q_k}{k}}+y^{\frac{q_k}{k}})^k\leq (x+y)^{q_k}.
\end{equation}

\begin{lemma}\label{elem lemma for higher energy}
For all $a\in[0,1]$ we have
\begin{equation}\label{goal ineq}
(a^{\frac{q_k}{k}}+(1-a)^{\frac{q_k}{k}})^k+2a^{\frac{q_k}{2}}(1-a)^{\frac{q_k}{2}}\leq 1.
\end{equation}

\end{lemma}
Observe that \eqref{Key elem ineq for higher energy} follows from \eqref{goal ineq} by taking $a=\frac{x}{x+y}$. 
A key ingredient in the proof of Lemma \ref{elem lemma for higher energy} is the following result established by Carlen, Frank, Ivanisvili and Lieb \cite[Proposition 3.1]{CFIL}.

\begin{proposition}\label{prop: CFIL}
For all $a\in[0,1]$ and $p\in (-\infty, 0]\cup[1,2]$
\begin{equation}\label{CFIL ineq}
(a^p+(1-a)^p)\left(1+\left(\frac{2a^{\frac{p}{2}}(1-a)^{\frac{p}{2}}}{a^{p}+(1-a)^p}\right)^{\frac{2}{p}}\right)^{p-1}\leq 1.  
\end{equation}
Moreover, the reverse inequality holds if $p\in[0, 1] \cup[2,\infty)$.
\end{proposition}

\begin{proof}[Proof of Lemma \ref{elem lemma for higher energy}]
We observe that \eqref{goal ineq} is equivalent to proving 
\begin{equation*}
1+\left(\frac{2^{\frac{1}{k}}a^{\frac{q_k}{2k}}(1-a)^{\frac{q_k}{2k}}}{a^{\frac{q_k}{k}}+(1-a)^{\frac{q_k}{k}}}\right)^k\leq \frac{1}{(a^{\frac{q_k}{k}}+(1-a)^{\frac{q_k}{k}})^k}.
\end{equation*}

Since $k<q_k=\log_{2}(2^k+2)<k+1$ for all $k\geq 2$, by taking $p=\frac{q_k}{k}$ in Proposition \ref{prop: CFIL} we obtain 
\begin{equation}\label{CFIL ineq consequence}
\left(1+\left(\frac{2a^{\frac{q_k}{2k}}(1-a)^{\frac{q_k}{2k}}}{a^{\frac{q_k}{k}}+(1-a)^\frac{q_k}{k}}\right)^{\frac{2k}{q_k}}\right)^{\frac{q_k}{k}-1}\leq \frac{1}{(a^{\frac{q_k}{k}}+(1-a)^{\frac{q_k}{k}})}.    
\end{equation}
Thus, it is enough to prove 
\begin{equation*}
1+\left(\frac{2^{\frac{1}{k}}a^{\frac{q_k}{2k}}(1-a)^{\frac{q_k}{2k}}}{a^{\frac{q_k}{k}}+(1-a)^{\frac{q_k}{k}}}\right)^k\leq \left(1+\left(\frac{2a^{\frac{q_k}{2k}}(1-a)^{\frac{q_k}{2k}}}{a^{\frac{q_k}{k}}+(1-a)^\frac{q_k}{k}}\right)^{\frac{2k}{q_k}}\right)^{q_k-k}.
\end{equation*}
Defining $\mu:=\frac{2a^{\frac{q_k}{2k}}(1-a)^{\frac{q_k}{2k}}}{a^{\frac{q_k}{k}}+(1-a)^\frac{q_k}{k}}$ (observe that $\mu\in[0,1]$ by AM-GM inequality), it is enough to prove 
$$
1+\frac{\mu^k}{2^{k-1}}\leq (1+\mu^{\frac{2k}{q_k}})^{q_k-k}
$$
for all $\mu\in[0,1]$. By letting $z:=\mu^{\frac{2k}{q_k}}$, we reduce the problem to proving
\begin{equation}
    \label{eq:concave_convex}
    1+\frac{z^{\frac{q_k}{2}}}{2^{k-1}}\leq (1+z)^{q_k-k}    
\end{equation}
for all $z\in[0,1]$. The equality holds at $z=0$ and $z=1$. 
Moreover, the left hand side of \eqref{eq:concave_convex} is convex in $z$ (as $2\leq k< q_k$),
and the right hand side is concave (as $k< q_k<k+1$). 
Therefore \eqref{eq:concave_convex} holds for all $z \in [0,1]$.

\end{proof}

\section{Proof of Theorem \ref{main thm 1}}

In this section  we show how to obtain Theorem \ref{main thm 1} from Lemma \ref{elem}, and then we  prove this lemma.
As before, we proceed via induction. Clearly, the result holds for $d=0$. Assume now  $d \geq 1$, and  the result has been established for $d-1$. Any set $A \subseteq \{0,1\}^d$ can be written as
$$ A = (A_0 \times \{0\}) \uplus (A_1 \times \{1\})$$
for some $A_0,A_1 \subseteq \{0,1\}^{d-1}$.

\ We have
\begin{align}
    \label{key formula}
    E_k(A)&= E_k(A_0)+E_k(A_1)\nonumber\\
    &\ \ \ + \sum_{i=1}^{k-1} {\binom{k}{i}}^2|\{ {(a_1,a_2,\dots,a_{2k})} \in A^i_0\times A^{k-i}_1 \times A^{i}_0 \times A^{k-i}_1\nonumber\\
&\qquad\qquad\qquad : a_1 +\dots+a_k = a_{k+1} +\dots + a_{2k} \}|\nonumber\\
&=E_k(A_0)+E_k(A_1)+ \sum_{i=1}^{k-1} {\binom{k}{i}}^2C_{i,k}.
\end{align}


Similarly to Proposition \ref{bounds for c i k 2}, we have
\begin{proposition}\label{bounds for c i k}
For all $1\leq i\leq k-1$ the following inequality holds
\begin{equation}\label{bounds}
C_{i,k}\leq |A_0|^{\frac{i}{k}p_k}|A_{1}|^{\frac{k-i}{k}p_k}.
\end{equation}
\end{proposition}

Observe that  Theorem \ref{main thm 1} follows from Proposition \ref{bounds for c i k}. 
Indeed, by \eqref{key formula}, Proposition \ref{bounds for c i k} and \eqref{elem ineq} we have
\begin{align*}
    E_{k}(A)&=E_k(A_0)+E_k(A_1)+ \sum_{i=1}^{k-1}{\binom{k}{i}}^2C_{i,k}\\
    &\leq E_k(A_0)+E_k(A_1)+ \sum_{i=1}^{k-1}{\binom{k}{i}}^2|A_0|^{\frac{i}{k}p_k}|A_{1}|^{\frac{k-i}{k}p_k}\\
    &\leq (|A_0|+|A_1|)^{p_k}\\
    &=|A|^{p_k}.
\end{align*}

\begin{proof}[Proof of Proposition \ref{bounds for c i k}]
We observe that
\begin{align*}
    C_{i,k}=\sum_{x\in\Z^d}|\chi_{A_0}*_{i-1}\chi_{A_0}*\chi_{A_1}*_{k-i-1}\chi_{A_1}(x)|^2,
\end{align*}
where, for compactly supported $f,g$, we define  $f*g(x):=\sum_{y\in\Z^d}f(y)g(x-y)$ and $*_{k}:=*(*_{k-1})$. Indeed, this follows from the fact that
$$\chi_{A_0}*_{i-1}\chi_{A_0}*\chi_{A_1}*_{k-i-1}*\chi_{A_1}(x)$$ counts the number of  $k$-tuples $(a_1,a_2,\dots,a_i,a_{i+1},\dots,a_k)\in A^i_0\times A^{k-i}_1$ such that $a_1+a_2+\dots+a_k=x$. Then, by Plancherel's theorem and H\"older's inequality we obtain
\begin{align*}
    C_{i,k}&=\sum_{x\in\Z^d}|\chi_{A_0}*_{i-1}\chi_{A_0}*\chi_{A_1}*_{k-i-1}\chi_{A_1}(x)|^2\\
    &=\int_{{\mathbb T}^d}|\Hat{\chi}_{A_0}(y)|^{2i}|\hat{\chi}_{A_1}(y)|^{2(k-i)}d m(y)\\
    &\leq \left(\int_{{\mathbb T}^d}|\Hat{\chi}_{A_0}(y)|^{2k}d m(y)\right)^{\frac{i}{k}}\left(\int_{{\mathbb T}^d}|\Hat{\chi}_{A_1}(y)|^{2k}d m(y)\right)^{\frac{k-i}{k}}\\
    &=\left( \sum_{x\in\Z^d}|\chi_{A_0}*_{k-1}\chi_{A_0}|^2\right)^{\frac{i}{k}}\left( \sum_{x\in\Z^d}|\chi_{A_1}*_{k-1}\chi_{A_1}|^2\right)^{\frac{k-i}{k}}\\
    &=E^{\frac{i}{k}}_k(A_0)E^{\frac{k-i}{k}}_k(A_1) \leq |A_0|^{\frac{ip_k}{k}}|A_1|^{\frac{(k-i)p_k}{k}},
\end{align*}
where $m$ is the Haar measure on $\mathbb{T}^d$ with $m(\mathbb{T}^d)=1$. 
\end{proof}

\begin{proof}[Proof of Lemma \ref{elem}] 
After re-scaling, we observe that to prove (\ref{elem ineq}) it is sufficient to show
\begin{equation}\label{key ineq}
\sum_{i=0}^{k}{\binom{k}{i}}^2x^{ip_k/k}  \leq (1+x)^{p_k}
\end{equation}
for all $1\leq x<\infty$. Moreover, after a change of variable, this is equivalent to proving that
\begin{equation}\label{key ineq 2}
g_k(y):=\sum_{i=0}^{k}{\binom{k}{i}}^2y^{i}\leq (1+y^{\alpha})^{\frac{k}{\alpha}}=:h_k(y)
\end{equation}
for all $1 \leq y \leq \infty$, where $\alpha:=\frac{k}{p_k}\in(1/2,1)$. Let $f(y):=\log h_{k}(y)-\log g_k(y)$. We need to show $f(y)\geq 0$ for all $y\geq 1$. Observe that 
$f(1)=0$. Moreover
$$
\lim_{y\to\infty} f(y)=\lim_{y\to\infty}\log\left(\frac{(\frac{1}{y^{\alpha}}+1)^{\frac{k}{\alpha}}}{\sum_{i=0}^{k}{\binom{k}{i}}^2y^{-i}}\right)=0,
$$
and, since
$$
\left(\frac{1}{y^{\alpha}}+1\right)^{\frac{k}{\alpha}}\geq 1+\frac{k}{\alpha y^{\alpha}}\ \ \text{and}\ \ \sum_{i=0}^{k}{\binom{k}{i}}^2y^{-i}=1+O(\frac{1}{y}),
$$
we have  $f(y)>0$ whenever $y$  is sufficiently large. Thus, it is sufficient to prove that $f'$ changes sign at most once in $(1,\infty)$. Observe that
\begin{align*}
yf'(y)&=\frac{ky^{\alpha}}{1+y^{\alpha}}-\frac{\sum_{i=0}^{k}{\binom{k}{i}}^2iy^i}{\sum_{i=0}^{k}{\binom{k}{i}}^2y^i} \\
&=\frac{y^{\alpha}\sum_{i=0}^{k}{\binom{k}{i}}^2y^i(k-i)-\sum_{i=0}^{k}{\binom{k}{i}}^2iy^i}{(1+y^{\alpha})\left(\sum_{i=0}^{k}{\binom{k}{i}}^2y^i\right)}.
\end{align*}

Thus, we need to prove that $y^{\alpha}\sum_{i=0}^{k}{\binom{k}{i}}^2y^i(k-i)-\sum_{i=0}^{k}{\binom{k}{i}}^2iy^i$ changes sign in $(1,\infty)$ at most once.
We define $$
\phi(y):=\log\left(y^{\alpha}\sum_{i=0}^{k}{\binom{k}{i}}^2y^i(k-i)\right)-\log\left(\sum_{i=0}^{k}{\binom{k}{i}}^2iy^i\right).
$$
We then have $\phi(1)=0$ and 
$$
\phi(y)=\alpha\log(y)+\log\left(\frac{n^2y^{n-1}+O(y^{n-2})}{ny^n+O(y^{n-1})}\right) \quad \text{as} \quad  y\to \infty.$$
Hence $\lim_{y\to\infty}\phi(y)=-\infty$\footnote{Here  we use the notation $V(y)=O(U(y))$ at $y_0$ to denote that an estimate of the form $|V(y)|\leq C |U(y)|$,   with some  constant  $C>0$, holds around $y_0$.}. It suffices to show that $\phi'$ changes sign (from + to -) at most once in $(1,\infty)$. Observe that
\begin{align*}
    \phi'(y)&=\frac{\alpha}{y}+\frac{\sum_{i=0}^{k}{\binom{k}{i}}^2y^{i-1}(k-i)i}{\sum_{i=0}^{k}{\binom{k}{i}}^2y^i(k-i)}-\frac{\sum_{i=0}^{k}{\binom{k}{i}}^2i^2y^{i-1}}{\sum_{i=0}^{k}{\binom{k}{i}}^2iy^i}\\
    &=\frac{\sum_{i=0}^{2k}C_iy^i}{y\left(\sum_{i=0}^{k}{\binom{k}{i}}^2y^i(k-i)\right)\left(\sum_{i=0}^{k}{\binom{k}{i}}^2iy^i\right)},
\end{align*}
where
\newcommand{\that}{\qquad \text{that}}
\begin{align*}
    C_i&:=\sum_{\substack{j+l=i\\0\leq j,l\leq k}}{\binom{k}{j}}^2{\binom{k}{l}}^2[\alpha(k-l)j+(k-l)lj-j^2(k-l)] \\  
    &=\sum_{\substack{j+l=i\\0\leq j,l\leq k}}{\binom{k}{j}}^2{\binom{k}{l}}^2 j(k-l)(\alpha +l-j) \,\,  \text{for all} \,\,  i,\,  0\leq i\leq 2k.
\end{align*}
Let $P(y):=\sum_{i=0}^{2k}C_iy^i$. We would like to show that $P(y)$ changes sign at most once from $+$ to $-$ in $(1,\infty)$. 
First, we claim $P(y)$ is a palindromic polynomial, i.e., $C_{i}=C_{2k-i}$ for all $i=0,\ldots, k$. Indeed, 
\begin{align*}
    &C_{2k-i} = \sum_{\substack{j+l=2k-i\\0\leq j,l\leq k}}{\binom{k}{j}}^2{\binom{k}{l}}^2 j(k-l)(\alpha +l-j)=  \\
    &\sum_{\substack{(k-j)+(k-l)=i\\0\leq j,l\leq k}}{\binom{k}{k-j} }^2{\binom{k}{k-l}}^2 (k-(k-j))(k-l)(\alpha +(k-j)-(k-l)). 
\end{align*}
If we denote $\tilde{l} =k-j$ and $\tilde{j}=k-l$, then we obtain 
\begin{align*}
   C_{2k-i} = \sum_{\substack{\tilde{l}+\tilde{j}=i\\0\leq \tilde{j},\tilde{l}\leq k}}{\binom{k}{\tilde{l}}}^2{\binom{k}{\tilde{j}}}^2 \tilde{j}(k-\tilde{l})(\alpha +\tilde{l}-\tilde{j}),
\end{align*}
which coincides with $C_{i}$. 
Since $P$ is the palindromic polynomial it follows that $y_{0}$ is its positive root if and only if  $P(1/y_{0})=0$. Therefore, to show that $P(y)$ changes sign from $+$ to $-$ at most once in $(1, \infty)$, it suffices to verify that $P(y)$ has at most two roots in $(0, \infty)$. By Descartes' rule of sign change $P(y)$ has at most two positive roots if there is at most two sign changes between consecutive (nonzero) coefficients $C_{i}$, $0\leq i \leq 2k$. Since $C_{i}=C_{2k-i}$ it suffices to show that there is at most one sign change between  consecutive (nonzero) coefficients, $C_{i}$ for $0\leq i \leq k$. 
Since $C_{0}=0$ we should consider coefficients $C_{i}$ with $1\leq i \leq k$. In the table below $C^{*}_{i} := \mathrm{sign}(C_{i})$, and $2\leq k \leq 10$. 
\begin{center}
\begin{tabular}{|c c c c c c c c c c c|} 
 \hline
 $k$ & $C_{1}^{*}$ & $C_{2}^{*}$ & $C_{3}^{*}$ & $C_{4}^{*}$ & $C_{5}^{*}$ & $C_{6}^{*}$ & $C_{7}^{*}$ & $C_{8}^{*}$ & $C_{9}^{*}$ & $C_{10}^{*}$ \\ [0.5ex] 
 \hline\hline
 2 & -1 & 1 &  &  &  &  &  &  &  &   \\ 
 \hline
 3 & -1 & 1 & 1 &  &  &  &  &  &  &   \\ 
 \hline
 4 & -1 & 1 & 1 & 1 &  &  &  &  &  &   \\ 
 \hline
 5 & -1 & 1 & 1 & 1 & 1 &  &  &  &  &   \\ 
  \hline
 6 & -1 & 1 & 1 & 1 & 1 & 1 &  &  &  &   \\ 
  \hline
 7 & -1 & -1 & 1 & 1 & 1 & 1 & 1 &  &  &   \\ 
  \hline
 8 & -1 & -1 & 1 & 1 & 1 & 1 & 1 & 1 &  &   \\ 
  \hline
 9 & -1 & -1 & -1 & 1 & 1 & 1 & 1 & 1 & 1 &   \\ 
  \hline
 10 & -1 & -1 & -1 & 1 & 1 & 1 & 1 & 1 & 1 & 1  \\ 
 \hline
\end{tabular}
\end{center}
\end{proof}
\begin{remark}\label{rem:conj}
It seems to us that Lemma \ref{elem} holds for all $k\geq 2$. 
We have verified at most one sign flip of the numbers $C_i$, $1\leq i \leq k$  on a computer for $k\le 100$.
It is an interesting question to verify that there is at most one sign flip in the sequence of $C_i$ for all $k$.
\end{remark}

\subsection*{Note added in proof} Motivated by Remark \ref{rem:conj}, Vjekoslav Kova\v{c} recently proved inequality \eqref{elem ineq} for all $k\geq 2$, see \cite{K}.

\begin{figure}
    \centering
\ref{fig:q_k}    \includegraphics[width=.74\linewidth, height=.44\linewidth]{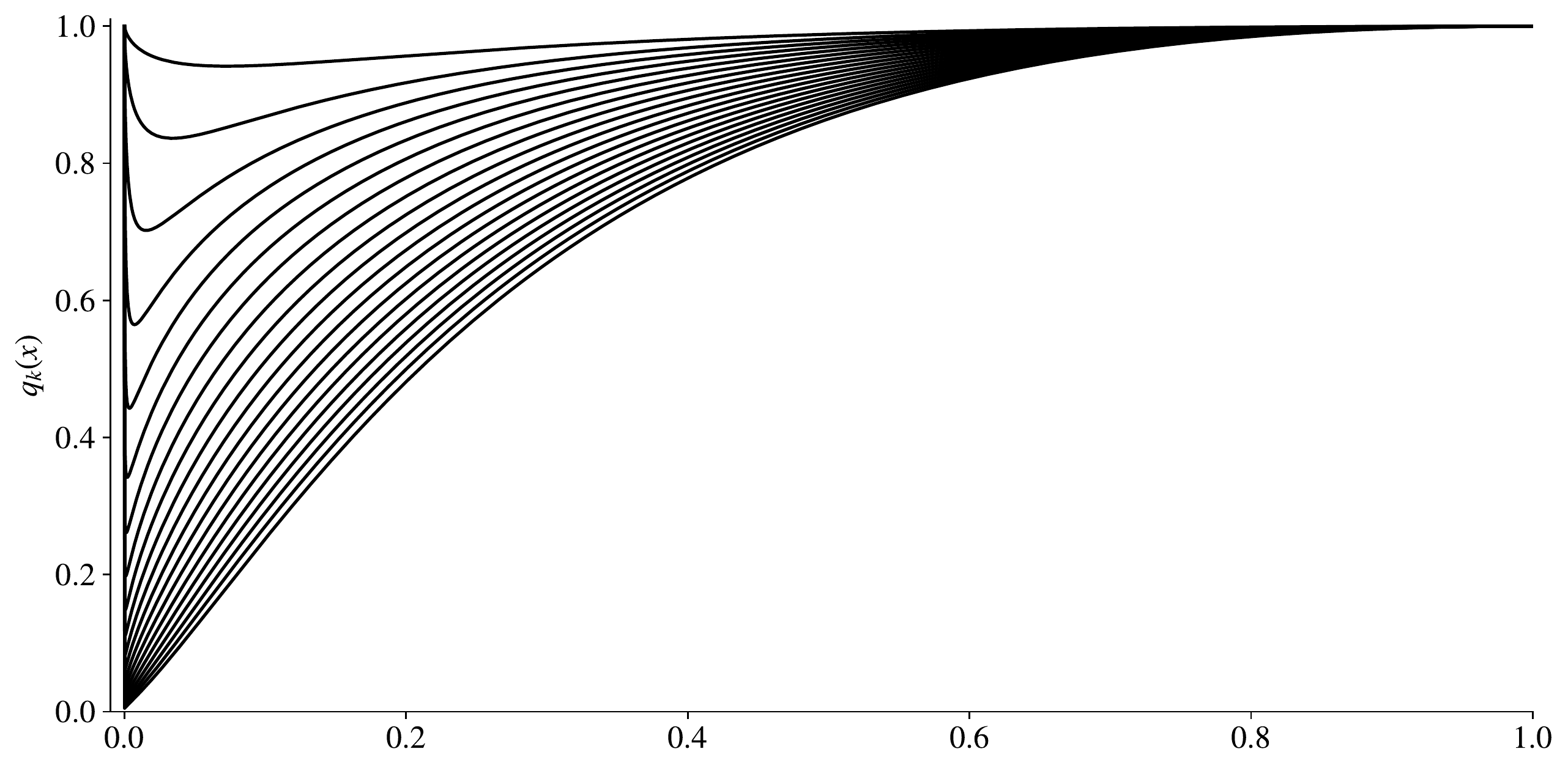}
    \caption{Graphs of $q_k(x)$  for $k\in\{2^{n};1\leq n< 20\}$. The picture suggests that $q_k(x)\leq 1$ for all $x\in[0,1]$. Lower graphs correspond to larger values of $k$.}
    \label{fig:q_k}
\end{figure}

\begin{figure}
    \centering
    \includegraphics[width=.74\linewidth, height=.44\linewidth]{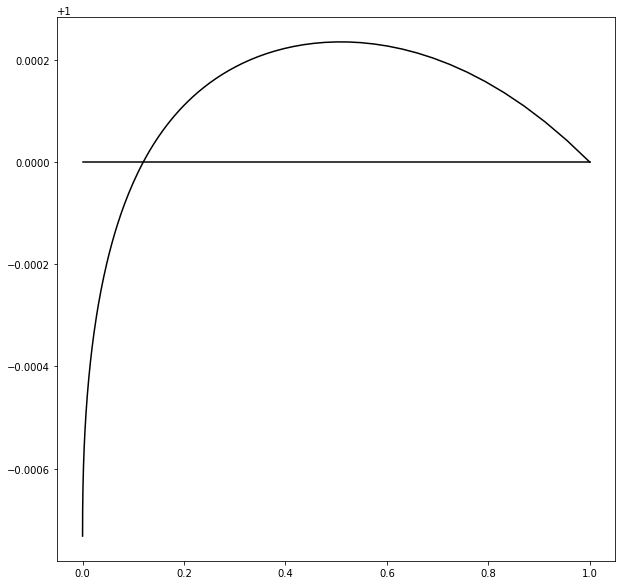}
    \caption{Graph of $\psi_3(x)$. We observe that $\psi_3(x)$ is  concave, and $\psi_3(x)$ intersects the line $y=1$ at only one point in $(0,1)$. }
    \label{fig:r_3}
\end{figure}

\begin{figure}
    \centering
    \includegraphics[width=.74\linewidth, height=.44\linewidth]{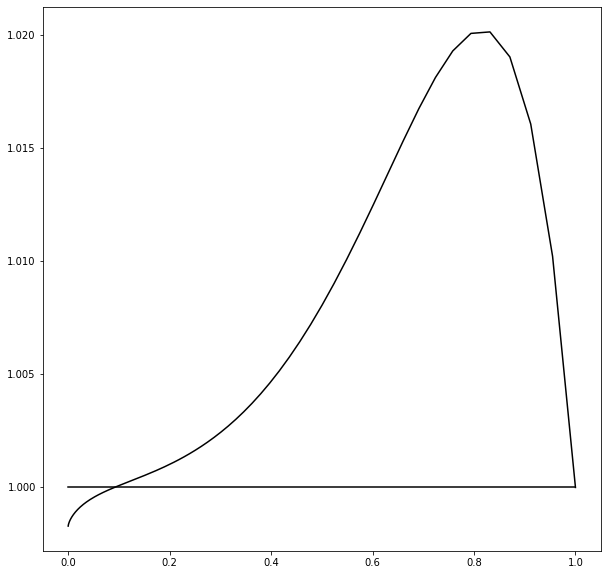}
    \caption{Graph of $\psi_7(x)$. We observe that $\psi_7(x)$ is not concave, however $\psi_7(x)$ still intersects the line $y=1$ at only one point in $(0,1)$. }
    \label{fig:r_7}
\end{figure}

\begin{remark}
To prove \eqref{key ineq} it suffices to show  
\begin{equation}\label{desired ineq}
\phi_k(x):=\frac{\sum_{i=0}^{k}{\binom{k}{i}}^2x^{p_k(k-i)/k}}{(1+x)^{p_k}}\leq 1
\end{equation}
for all $x\in[0,1]$.
The inequality (\ref{desired ineq})  can be easily verified around $x=0$. One can also verify it around $x=1$. Therefore, to obtain the desired inequality in the whole interval $[0,1]$ it would be enough to prove that each $\phi_k$ has only one critical point in $(0,1)$. We observe that $x$ is a critical point of $\phi_k$ if and only if 
\begin{equation*}
(1+x)^{p_k+1} \phi_k'(x) = 
\sum_{i=0}^{k}{\binom{k}{i}}^2\left[\frac{p_k(k-i)}{k}x^{p_k\frac{k-i}{k}-1}(1+x)-p_k x^{p_k\frac{k-i}{k}}\right]=0,
\end{equation*}
or, equivalently
\begin{equation*}
\psi_{k}(x):=\sum_{i=0}^{k-1}{\binom{k}{i}}^2\left[\frac{k-i}{k}x^{p_k\frac{k-i}{k}-1}-\frac{i}{k}x^{p_k\frac{k-i}{k}}\right]=1.
\end{equation*}
Therefore, as $\psi_k(0)=0$ and $\psi_k(1)=1$, in order to establish the desired inequality, i.e., $\phi_k(x)\leq1$ for all $x\in(0,1)$, it would be enough to prove that $\psi_k(x)$ is concave.  For small values of $k$, one can establish the concavity of $\psi_k$; in particular, this is the approach of Kane--Tao \cite{KT} for $k=2$. Figure \ref{fig:r_3} illustrates that $\psi_k$ is concave for $k=3$. Unfortunately, this is no longer the case if $k$ is large; e.g., Figure \ref{fig:r_7} illustrates the non-concavity of $\psi_k$  for $k$ as small as $7$ already. Another approach to prove Lemma \ref{elem} would be to show $\phi_{k+1}(x)\leq \phi_k(x)$ which numerically  seems correct. 
\end{remark}


\smallskip

\section{Proofs of Propositions \ref{Thm: 0,1,2} and  \ref{Thm larger cubes}}

\newcommand{\DR}{\operatorname{DE}}

The proof of Kane--Tao  \cite{KT} of the $\{0,1\}$-analogue, 
as well as the proofs of Theorems \ref{main thm 2} and \ref{main thm 1}   are based on the following two steps:

\begin{itemize}
    \item \emph{Guessing} the extremizer to the inequality (which, in those cases, happened to be the entire set).
    \item Showing an inductive bound that allowed us to see that the extremizer candidate is indeed the extremizer.
\end{itemize}

In the $\{0,1,2\}^n$ or more general cases the entire set is not generally the extremizer, 
and finding the extremizer becomes a key step of the proof:

\begin{itemize}
    \item We first construct an auxiliary problem that \emph{inducts}, or, in this case \emph{tensorizes} 
    essentially by construction. 
    Solving this problem is essentially equivalent to  \emph{guessing} the extremizers in the previous problems.
    \item We then show that the solution to this auxiliary problem gives rise to sharp (almost) extremizers of the original problem. 
    This step is new, and necessary due to the fact that the extremizing sets are in 
    general far from being product sets. 
\end{itemize}

\subsection{The auxiliary (discrete restriction) problem}

For each specific instance of interest (in our case $\{0,1,2\}$) 
the auxiliary problem will then reduce to solving a finite-dimensional optimization problem
closely related to the inequalities studied in the previous sections.
The way to define these problems will be by defining auxiliary quantities frequently appearing in the  discrete restriction theory.

\begin{definition}[Discrete extension constants]
    Given positive integers $k,d$, and a finite subset $A\subset \mathbb Z^d$, we define:

    \begin{itemize}
        \item The discrete extension constant $\DR_{l^q\to L^{2k}}(A)$ as the smallest constant such that, 
        for any function $f:A\to \mathbb R$ it holds that
        \begin{equation}
        \label{eq:DR_def}
            \left |
                \sum_{\substack{x_1, \dots, x_k \in A \\ y_1, \dots, y_k \in A \\ \sum x_i = \sum y_i}} 
                    \hspace{-1em} 
                     f(x_1)\cdot \dots \cdot  f(x_k) \cdot  f(y_1)\dots f(y_k)
            \right |^{\frac 1 {2k}} \hspace{-.8em} 
            \le 
            \DR_{l^q \to L^{2k}}(A) \|f\|_{l^q(A)}.
        \end{equation}
        \item The restricted discrete extension constant $\tilde {\DR}_{l^q\to L^{2k}}(A)$, 
        which is the best possible  constant so that \eqref{eq:DR_def} holds for all functions $f:A\to \{0,1\}$.
    \end{itemize}
\end{definition}

The quantities $\DR$, $\tilde \DR$ 
have essentially the same value (Lemma \ref{lem:comparison}), 
but $\DR$ is much easier to work with (Lemma \ref{lem:Tensorization}).
Moreover, understanding for which $q$ we have $\tilde \DR_{l^q \to L^{2k}}(\{0,1,2\}^d)\leq 1$ is  essentially equivalent to proving Proposition \ref{Thm: 0,1,2}.

\begin{lemma}\label{lem:convolution_FT}
    Let $A$ be a finite subset of $\mathbb Z^d$. Let $1\leq p = \frac{2k}{q}$, and $C>0$. 
    The following statements are equivalent:

    \begin{enumerate}
        \item For all subsets $B\subset A$, it holds that
         \begin{equation*}
             E_k(B) \le C^{2k}|B|^p.
         \end{equation*}
        \item 
        \begin{equation*}
            \tilde \DR_{l^q \to L^{2k}}(A) \le C.
        \end{equation*}
    \end{enumerate}
\end{lemma}

\begin{proof}
    Set $f$ in the definition of $\tilde \DR$ to be equal to $\chi_B$ for $B$ as in  part (1).
\end{proof}

The constant $\DR$ is called the discrete extension constant because it is, indeed, the operator norm of  an extension operator.

\begin{lemma}[Fourier transform]
\label{lem:fourier_convolution}
    Let $A$ be a finite subset of $\mathbb Z^d$. 
    Then $\DR_{l^q\to L^{2k}}(A)$ is the operator norm of the extension operator\footnote{Here we denote by $\mathcal F\{f\}$ the Fourier transform of $f$, i.e., $\mathcal F\{f\}(z)=\sum_{k\in \Z^d}f(k)z^k$.} $\mathcal E(f)=\mathcal F\{f\}$ 
    from $l^q(A)\subseteq l^q(\mathbb Z^d)$ to $L^{2k}(\mathbb T^d)$.
\end{lemma}

\begin{proof}
    By definition, $\DR_{l^q\to L^{2k}}(A)$ is the best constant such that, 
    for any function $f:\mathbb Z^d \to \mathbb R$ supported on $A$, it holds that:
    \begin{equation*}
        \|
            \underbrace{f\ast f\ast f \dots \ast f}_{k \text{ times}}
        \|_{l^2(\mathbb Z^d)}^{1/k}
        \le 
        \DR_{l^q\to L^{2k}}(A) \|f\|_{l^{q}(\mathbb Z^d)}.
    \end{equation*}
    At the same time, by Plancherel's theorem and the product-convolution rule
    \begin{align*}
        \|
            f^{*k}
        \|_{l^2(\mathbb Z^d)}
        = 
        \|
       \mathcal F \{ {f^{*k}}\}
        \|_{L^2(\mathbb T^d)}
        =
        \|
        \mathcal F \{ f\}^ k
        \|_{L^2(\mathbb T^d)}
        =
        \|
        \mathcal F \{ f\}
        \|_{L^{2k}(\mathbb T^d)}^{k}.
    \end{align*}

\end{proof}

\begin{remark}
    Lemma \ref{lem:fourier_convolution} above shows that the constants 
    $\tilde \DR_{l^q \to L^{2k}}(A), \DR_{l^q \to L^{2k}}(A)$ 
    make sense for arbitrary $2k \in \mathbb R$, and not just even integers. 
\end{remark}

The following  lemma is essentially \cite[Proposition 3.3]{dec}. For  completeness of the argument we include the proof here.

\begin{lemma}
    [Tensorization Lemma]
    \label{lem:Tensorization}
    Let $q\le 2k$. Then for 
        $A\subseteq \mathbb Z^{d_1}$, 
        $B\subseteq \mathbb Z^{d_2}$, 
        $A\times B \subseteq \mathbb Z^{d_1}\times \mathbb Z^{d_2}$ 
    we have

    \begin{equation*}
        \DR_{l^q \to L^{2k}} (A\times B) = \DR_{l^q \to L^{2k}} (A) \DR_{l^q \to L^{2k}}(B).
    \end{equation*}
\end{lemma}

\begin{proof}
    The ``$\ge$'' inequality follows by testing the left hand side operator 
    with the tensor product of (almost) extremizers to the right hand side. 

    For the opposite direction, let 
        $f:A\times B \to \mathbb C$
        , and 
        $\hat f: \mathbb T^{d_1} \times  \mathbb T^{d_2} \to \mathbb C$ 
    be its Fourier transform. 
    Let $\mathcal F_1, \mathcal F_{2}$ be the Fourier transforms on $\mathbb Z^{d_1}$ and $\mathbb Z^{d_2}$. 
    The goal is to estimate
    \begin{equation*}
            \|
                \|
                \mathcal F_{2}\{\mathcal F_1 f\}(x_1,x_2)
                \|_{L^{2k}(x_2\in \mathbb T^{d_2})}
            \|_{L^{2k}(x_1\in \mathbb T^{d_1})}.
    \end{equation*}
    Fixing $x_2$, we apply the $\DR$ inequality
    \begin{equation*}
        \|
            \mathcal F_{2}\{\mathcal F_1 f\}(x_1,x_2)
        \|_{L^{2k}(x_2\in \mathbb T^{d_2})} 
        \le  
        \DR_{l^q \to L^{2k}} (B) 
        \|
            \mathcal F_1 f(x_1,b)
        \|_{l^q(b)}.
    \end{equation*}
     Now, using the hypothesis that $2k\ge q$, we can reverse the norms
    \begin{equation*}
        \|
            \|
                \mathcal F_1 f(x_1,b)
            \|_{l^q(b\in  B)}
        \|_{L^{2k}(x_1\in \mathbb T^{d_1})}
        \le
        \|
            \|
                \mathcal F_1 f(x_1,b)
            \|_{L^{2k}(x_1\in \mathbb T^{d_1})}
        \|_{l^q(b\in  B)}.
    \end{equation*}
    Now the $\DR$ inequality can be applied again to $\|
                \mathcal F_1 f(x_1,b)
            \|_{L^{2k}(x_1\in \mathbb T^{d_1})}$. Joining it all together
    \begin{equation*}
    \begin{aligned}
        \|
            \mathcal F_{2}(\mathcal F_1 1)&(x_1,x_2)
        \|_{L^{2k}(x_1\in \mathbb T^{d_1}))
            L^{2k}(x_2\in \mathbb T^{d_2}))}
        \le\\&\le
        \DR_{l^q \to L^{2k}} (B)
        \DR_{l^q \to L^{2k}} (A)
        \|f(a,b)\|_{l^{q}(a\in A)l^q(b\in B)}.
    \end{aligned}
    \end{equation*}
\end{proof}

\subsection{Relating the Discrete extension problem and the original problem}

In this section we show that the discrete extension constants 
    $\tilde \DR_{l^q \to L^{2k}} (A^d)$ and $\DR_{l^q \to L^{2k}} (A^d)$ 
grow similarly as $d$ goes to infinity. This will allow us to compute 
the assymptotic behavior of $\DR$
in order to find the (much harder) assymptotics for $\tilde \DR$. The next lemma is inspired by Bourgain's logarithmic pigeonhole principle (see \cite{Tao}).

\begin{lemma}
    [Comparison Lemma]
    \label{lem:comparison}
    For all $q\ge 1$, $k\ge \frac 1 2$, $A\subseteq \mathbb Z^d$ it holds that
    \begin{equation*}
        \tilde \DR_{l^q \to L^{2k}} (A) 
        \le 
        \DR_{l^q \to L^{2k}} (A)  
        \le 
        (2+\log |A|)\tilde \DR_{l^q \to L^{2k}} (A).
    \end{equation*}
\end{lemma}

\begin{proof}
    The first inequality follows by the fact that $\DR$ is a maximum over a larger class of functions. 
    For the second one, let $f:A \to \mathbb R$. Without loss of generality assume $\|f\|_{l^\infty(A)} = 1$, and that $f$ is nonnegative. We can decompose $f$ as a sum
    \begin{equation*}
        f(x) = \sum_{\substack{i\geq 1 \\ 2^i\leq |A| }}  2^{-i} \epsilon_i(x) + f_0 (x)
    \end{equation*}
    with the property that  $\epsilon_i: A\to \{0,1\}$, and $0\leq  f_0(x) \le |A|^{-1}$. 
    The value of $\epsilon_i(x)$ is the $i-$th digit of the boolean expantion of $f(x)$. 
    Moreover, $\|f_0\|_1 \le 1$.
     There are, moreover at most $(\log |A|+1)$ different $\epsilon_i$. By the triangle inequality, we have 
    \begin{equation*}
        \|\hat f\|_{L^{2k}(\mathbb T^d)} \le  
       \sum_{\substack{i\geq 1 \\ 2^i\leq |A| }}  2^{-i} \|\hat \epsilon_i\|_{L^{2k}(\mathbb T^d)} + \|\hat f_0 \|_{L^{2k}(\mathbb T^d)}.
    \end{equation*}
    We bound the sum  by the maximum element in the sum (times the number of elements),
    and the  term  $\|\hat f_0 \|_{L^{2k}(\mathbb T^d)}$ by $1$,  to obtain
    \begin{equation*}
        \|\hat f\|_{L^{2k}(\mathbb T^d)} \le    (1+\log(|A|)) \max_{i\geq 1} 2^{-i} \|\hat \epsilon_i\|_{L^{2k}(\mathbb T^d)} + 1.
    \end{equation*}
    Now, by applying the $\tilde \DR$ bounds on $\hat \epsilon_i$ we get
    \begin{equation*}
        \|\hat f\|_{L^{2k}(\mathbb T^d)} \le    (1+\log(|A|)) \tilde \DR_{l^q \to L^{2k}} (A)    \max_{i\geq 1} 2^{-i} \| \epsilon_i\|_{L^{q}(A)} + 1.
    \end{equation*}
    By construction $2^{-i} \| \epsilon_i\|_{L^{q}(A)} \le \|f\|_{L^q(A)}$. 
    By checking against a singleton, $\tilde \DR$ is always at least $1$, and $ \|f\|_{l^q(A)} \ge  \|f\|_{l^\infty(A)}= 1$. Combining all this, we obtain
    \begin{equation*}
        \|\hat f\|_{L^{2k}(\mathbb T^d)} \le    (2+\log(|A|)) \tilde \DR_{l^q \to L^{2k}} (A)\|f\|_{l^q(A)}.
    \end{equation*}
\end{proof}

\begin{remark}
    The exponent of the $\log$ in Lemma \ref{lem:comparison} is probably not sharp (see, for example, the gains in the log-power in \cite[Theorem 1.1]{reversing} or \cite[Lemma 2.4]{Mudgal}). Finding the sharp exponent is not necessary for our purposes. We thank A. Mudgal for this remark.
\end{remark}

The results from this section yield the relationship between Proposition \ref{Thm: 0,1,2} and the discrete extension constant, as follows.

\begin{proposition}\label{thm:DR_implies_energy}
    Let $A$ be a finite subset of $\mathbb Z$. Let $1\leq p = \frac{2k}{q}$, and $C>0$. The following are equivalent:
    \begin{enumerate}
        \item An inequality of the form 
        \begin{equation*}
        E_k(X)\le C |X|^p
        \end{equation*}
        holds for all  $X\subseteq A^d$, $d\geq 0$.
        \item An inequality of the form 
        \begin{equation*}
        E_k(X)\le |X|^p.
        \end{equation*}
        holds for all  $X\subseteq A^d$,  $d\geq 0$.
        \item $\DR_{l^q \to L^{2k}} (A)\le 1$.
    \end{enumerate}
\end{proposition}

\begin{proof}
    Clearly, (3) $\Rightarrow$ (2) (by Lemma \ref{lem:Tensorization} and Lemma \ref{lem:convolution_FT}), and (2) $\Rightarrow$ (1). We  show that (1) $\Rightarrow$ (3). By Lemmas \ref{lem:comparison} and \ref{lem:Tensorization} we have 
    \begin{equation}
    \label{eq:comparison_in_thm}
        \tilde\DR_{l^q \to L^{2k}} (A^d)\le \DR_{l^q \to L^{2k}} (A)^d \le (2+d\log |A|)\tilde\DR_{l^q \to L^{2k}} (A^d).
    \end{equation}
    Observe that by Lemma \ref{lem:convolution_FT}, (1) is equivalent  to
    \begin{equation}
    \label{eq:limsup}
        \sup_{d}\tilde\DR_{l^q \to L^{2k}} (A^d) <\infty.
    \end{equation}
    By equation \eqref{eq:comparison_in_thm}, equation \eqref{eq:limsup} is equivalent to
    \begin{equation*}
        \DR_{l^q \to L^{2k}} (A)^d \le 1
    \end{equation*}
    and the result follows.
\end{proof}

\begin{remark}
    The proof of Theorem \ref{thm:DR_implies_energy} extends to any finite subset $A$ of an abelian group $G$ without any significant changes, 
    using that the group generated by $A$ inside of $G$ is locally compact and abelian with the discrete topology.
\end{remark}

\subsection{Concluding the proofs of  Propositions \ref{Thm: 0,1,2} and \ref{Thm larger cubes}}

\begin{proof}[Proof of Proposition \ref{Thm: 0,1,2}]

    Applying Proposition \ref{thm:DR_implies_energy} with $A= \{0,1,2\}\subseteq \mathbb Z$ and $k=2$ shows that $t_2$ is equal to the smallest $p$ such that
    \begin{align*}
        \frac{x^p+y^p+z^p+4(x^{p/2}y^{p/2}+x^{p/2}z^{p/2}+y^{p/2}z^{p/2})+4x^{p/2}y^{p/4}z^{p/4}}{(x+y+z)^p}\leq 1,
    \end{align*}
    for all $x,y,z\geq 0$. In particular, taking $x=1$ and $y=z=1/2$ we obtain
    \begin{align*}
        t_2&\geq \inf\{p\in [2,3]\colon \ 4^p-2^p-12(2^{p/2})-6\geq 0\}\\
        &=\inf\{2\log_{2}w\colon w\in[2,2\sqrt{2}],\ w^4-w^2-12w-6\geq 0\}\\
        &\geq 2\log_2(2.5664)> \log_319 =\frac{\log E(\{0,1,2\}^d)}{\log |\{0,1,2\}^d|}.
    \end{align*}

    \end{proof}

    \begin{proof}[Proof of Proposition \ref{Thm larger cubes}]
    The upper bound is trivial, so we focus our attention on the lower bounds.
    Consider the case $n=2m-1$. We prove that $E(\{0,1,\dots,2m-1\})= \frac{16m^3+2m}{3}$. We start observing that for any $a\in\{0,1,\dots,2m-1\}$ the 4-tuple $(a,a,a,a)$ is a solution. Moreover, for all $a,b\in\{0,1,\dots,2m\}$ we have that $(a,b,a,b),(a,b,b,a),(b,a,b,a)\ \text{and}\ (b,a,a,b)$ are also solutions. This gives a total of
    $2m+4{\binom{2m}2}$ trivial solutions.
    
    Then, we observe that the $m$ couples
    $(0,2m-1),(2,2m-3),(3,2m-4),\dots,(m-1,m)$ add up to $2m-1$, this gives 
    $8{\binom{m} 2}$ nontrivial solutions.
    Similarly, the couples adding up to $2m-2$ and $2m$ give $8{\binom{m-1}2}+4{\binom{m-1}1}$ solutions. More generally, we have that considering the couples adding $k$ or $4m-2-k$ we obtain $8{\binom{\lceil k/2\rceil} 2}$ non-trivial solutions if $k$ is odd and   $8{ \binom{k/2}2}+4{\binom{k/2}1}$ if $k$ is even. Therefore
    \begin{align*}
    &E_2(\{0,1,\dots,2m-1\})\\
    &=2m+4{\binom{2m}2}
    +8{\binom{m} 2}
    +4\left(8\sum_{k=2}^{m-1}{\binom{k} 2}\right)
    +2\left(4\sum_{k=1}^{m-1}k\right)\\
    &=\frac{16m^3+2m}{3}.
    \end{align*}
    The case $n=2m$ follows similarly.
    \end{proof}


\smallskip





\section*{Acknowledgments} 
We are grateful to the anonymous referees and to the editors of Discrete Analysis for thoughtful comments and corrections which improved the exposition of the paper. We are also thankful to Terence Tao for helpful discussions. 

This work was initiated at the Hausdorff Research Institute for Mathematics, during
the trimester program ``Harmonic Analysis and Analytic Number Theory''; we are grateful to the institute and the organizers of the program.

\bibliographystyle{amsplain}


\begin{dajauthors}
\begin{authorinfo}
[J. de Dios]
Jaume de Dios Pont\\
Department of  Mathematics,  University  of  California  Los  Angeles\\
Portola Plaza 520, Los  Angeles,
CA 90095, USA\\
 jdedios\imageat{}math.ucla\imagedot{}edu 
\end{authorinfo}

\begin{authorinfo}
[R. Greenfeld]
Rachel Greenfeld\\
Institute for Advanced Study\\
Princeton, NJ 08540, USA\\
 greenfeld.math\imageat{}gmail\imagedot{}com
\end{authorinfo}

	\begin{authorinfo}
[P. Ivanisvili]
Paata Ivanisvili\\
Department of Mathematics, University of California Irvine\\
Rowland Hall 510C, Irvine, CA 92697, USA\\
 pivanisv\imageat{}uci\imagedot{}edu
\end{authorinfo}
	
		\begin{authorinfo}
[J. Madrid]
Jos\'e Madrid\\
Department of  Mathematics,  University  of  California  Los  Angeles\\
Portola Plaza 520, Los  Angeles,
CA 90095, USA\\
 jmadrid\imageat{}math.ucla\imagedot{}edu
\end{authorinfo}
	
\end{dajauthors}

\end{document}